\documentclass[12pt]{amsart}

\usepackage[margin=1in]{geometry}

\usepackage{amssymb,amsmath,amsthm,mathrsfs,enumerate,multicol,esint}
\usepackage{color}
\usepackage{makecell}
\usepackage{tikz}
\usetikzlibrary{arrows,calc}
\usepackage{hyperref}

\allowdisplaybreaks

\title[$L^2$ boundedness of Hermite pseudo-multipliers]{On the $L^2$ boundedness of  pseudo-multipliers for Hermite expansions}

\author{Fu Ken Ly}
\address{School of Mathematics and Statistics, The Learning Hub, The University of Sydney,
NSW 2006, Australia.}
\email{ken.ly@sydney.edu.au}


\newcommand{\RR}{\mathbb{R}} 
\newcommand{\NN}{\mathbb{N}} 
\newcommand{\ZZ}{\mathbb{Z}} 
\newcommand{\sz}{\mathscr{S}} 

\newcommand{\N}{\mathcal{N}}
\newcommand{\K}{\mathcal{K}}

\newcommand{\f}{\frac}
\newcommand{\lesi}{\lesssim}

\newcommand{\LL}{\mathcal{L}}
\newcommand{\QQ}{\mathbb{Q}} 
\newcommand{\ip}[1]{\langle #1 \rangle} 
\newcommand{\bip}[1]{\big\langle #1 \big\rangle} 
\newcommand{\Bip}[1]{\Big\langle #1 \Big\rangle} 

\newcommand{\wt}[1]{\widetilde{#1}} 
\newcommand{\diff}{\triangle}

\newcommand{\SM}{\mathcal{S}} 
\newcommand{\floor}[1]{\lfloor #1 \rfloor}

\newcommand{\A}[1]{A^{(#1)}}
\newcommand{\vph}{\varphi}
\newcommand{\dd}{\delta} 
\newcommand{\rr}{\rho} 

\newcommand{\D}{\partial}
\newcommand{\half}{\frac{1}{2}}

\newcommand{\I}{\mathcal{I}}
\newcommand{\ep}{\epsilon}

\newcommand{\KK}{K}

\newcommand{\om}{\omega}
\newcommand{\HH}{\mathcal{H}}
\newcommand{\OU}{\mathbb{L}} 

\newcommand{\U}{\mathcal{U}}

\DeclareMathOperator{\supp}{supp\,} 

\theoremstyle{plain}
\newtheorem{Theorem}{Theorem}[section]
\newtheorem{thm}{Theorem}
\newtheorem{Lemma}[Theorem]{Lemma}
\newtheorem{Proposition}[Theorem]{Proposition}
\newtheorem{Corollary}[Theorem]{Corollary}

\newtheorem{Definition}[Theorem]{Definition}
\theoremstyle{definition}
\newtheorem{Remark}[Theorem]{Remark}
\theoremstyle{remark}

\theoremstyle{example}

\numberwithin{equation}{section}

\def\barint{\kern4pt
\raise3.4pt\hbox{\vrule height.8pt width5pt}%
\kern-9pt 
\int}

\def\XXint#1#2#3{{\setbox0=\hbox{$#1{#2#3}{\int}$}
     \vcenter{\hbox{$#2#3$}}\kern-.5\wd0}}

\begin{document}

\begin{abstract}
We give various conditions for Hermite pseudo-multipliers to be bounded on $L^2(\mathbb{R}^n)$. As a by-product we also give  new results for pseudo-multipliers in the Gaussian measure setting. One of our key tools is a new integration by-parts formula for Hermite expansions.
\end{abstract}

\maketitle

\section{Introduction}\label{sec: intro}

The classical theory of pseudo-differential operators on $\RR^n$ is well established; such operators are defined by 
\begin{align}\label{eq:pdo}
\sigma(x,D) f(x)=\int \sigma(x,\xi)\,\widehat{f}(\xi) e^{2\pi i x\cdot \xi}\,d\xi
\end{align}
where $\sigma:\RR^n\times\RR^n\to \mathbb{C}$ is a bounded and measurable function. When $\sigma$ is independent of $\xi$, the operator is called a Fourier multiplier, and through the Plancherel theorem, is bounded on $L^2(\RR^n)$ if and only if the symbol is a bounded function. Sufficient conditions for pseudodifferential operators were given in the celebrated Calderon--Vaillancourt theorem, which states that \eqref{eq:pdo} extends to a bounded operator on $L^2(\RR^n)$ provided
\begin{align}\label{eq:CVpdo}
|\D_x^\nu\D_\xi^\kappa \sigma(x,\xi)|\le C_{\nu,\kappa} \ip{\xi}^{-\rr|\kappa|+\dd|\nu|}, \qquad |\nu|, |\kappa|\le \floor{n/2}+1
\end{align}
for $0\le \dd\le \rr\le 1$ and $\dd\ne 1$ (see \cite{CV}). Note that here and throughout the rest of the article we adopt the notation $\ip{\xi}=1+|\xi|$.

In this article we consider the $L^2$ boundedness of pseudo-differential type operators associated with certain orthogonal expansions derived from the well-known Hermite polynomials $\{H_k\}_{k\in \NN_0}$ which are given by  
$$H_k(t)=(-1)^ke^{t^2}\D_t^k(e^{-t^2}), \qquad t\in\RR,$$
and $\NN_0:=\NN\cup\{0\}$. 
We shall consider two kinds of expansions; the first type concerns so-called Hermite functions and the second are the normalised form of the Hermite polynomials.

\subsection{Pseudo-multipliers associated with Hermite functions}
For $n=1$ the Hermite function of degree $k\in\NN_0$ is
\begin{align*}
	h_k(t)=(2^k k! \sqrt{\pi})^{-1/2}H_k(t)e^{-t^2/2} \qquad \forall t\in\RR,
\end{align*}
while for $n\ge 2$ the $n$-dimensional Hermite functions $h_\xi$ are defined over the multi-indices $\xi\in \NN_0^n$ by 
\begin{align*}
	h_\xi(x)=\prod_{j=1}^n h_{\xi_j}(x_j)\qquad \forall x\in\RR^n.
\end{align*}
These functions are eigenfunctions of the operator $\LL=-\Delta+|x|^2$ in the sense that $\LL(h_\xi)=~(2|\xi|+~n)h_\xi$ for every $\xi\in\NN_0^n$; 
furthermore, they form an orthonormal basis for $L^2(\RR^n)$.

Given a bounded function $\sigma: \RR^n\times \NN_0^n\to\mathbb{C}$, we define the Hermite \emph{pseudo}-multiplier  $\sigma(\cdot,\LL)$ by
\begin{align}\label{eq:hpdo}
\sigma(x,\LL) f(x)
=\sum_{\xi \in \NN_0^n}\sigma(x,\xi)\ip{f,h_\xi} h_\xi(x),
\end{align}
where $\LL=-\Delta+|x|^2$ is the Harmonic Oscillator. Pseudo-multipliers were first introduced in \cite{Epp} and further investigated in \cite{BT,Bui,CR,LN21,LN22}. 

In the case that $\sigma$ is independent of $x$ (that is, $\sigma(x,\LL)=\sigma(\LL)$ is a Hermite multiplier) it is well known using  Parseval's identity that $\sigma\in L^\infty(\NN_0^n)$ is a necessary and sufficient condition for the $L^2(\RR^n)$ boundedness of the multiplier $\sigma(\LL)$.
However, unlike the multiplier case, it is not clear whether $\sigma(\cdot,\LL)$ is automatically bounded on $L^2(\RR^n)$ if $\sigma\in L^\infty(\RR^n\times \NN_0^n)$.
Indeed the following question was raised in \cite{BT}:
\begin{quote}
{\it What conditions on $\sigma(x,\xi)$  ensure that the pseudo-multiplier $\sigma(\cdot,\LL)$ extends to a  bounded operator on $L^2(\RR^n)$?}
\end{quote}
\noindent In particular,  the possibility of a Calder\'on--Vaillencourt type theorem for Hermite expansions was suggested.

One answer  was given in \cite{LN21}; the result there was modelled on the classical  `H\"ormander-type' condition \eqref{eq:CVpdo}. Indeed, in \cite{LN21} the following H\"ormander-type class was introduced. 
\begin{Definition}
Let $m\in\RR$, $\rr,\dd\ge0$ and $\N,\K\in \NN_0\cup\{\infty\}$. The symbol $\sigma:\RR^n\times \NN_0^n\to\mathbb{C}$ belongs to $\SM^{m,\K,\N}_{\rr,\dd}$ if $\sigma(\cdot,\xi)\in C^\N(\RR^n)$ for each $\xi \in \NN_0^n$ and there exists $C_{\nu,\kappa}>0$ such that for each $(x,\xi)\in\RR^n\times\NN_0^n$ we have 
\begin{align*}
	| \D^\nu_x \diff^\kappa_\xi \sigma(x,\xi)| \le C_{\nu,\kappa} \,\ip{\xi}^{\f{m}{2}-\rr|\kappa|+ \f{\dd}{2} |\nu|}  
	\qquad |\kappa|\le \K,\; |\nu|\le\N
\end{align*}
for $\nu,\kappa\in\NN_0^n$ satisfying $0\le |\nu|\le\N$ and $0\le|\kappa|\le \K$. 
\end{Definition}
When $\N=\K=\infty$ we just write $\SM^{m,\infty,\infty}_{\rr,\dd}=\SM^m_{\rr,\dd}$. 
 Note that here and in the sequel, if
$g$ is a function defined over $\NN_0^n$ then we set
$ \diff_i g(\xi) := g(\xi+e_i)-g(\xi)$,
and 
$ \diff_i^\ell g := \diff(\diff_i^{\ell-1}g)$  for $\ell\ge 2$. 

Then in \cite{LN21} the following was obtained. 

\begin{thm}[\cite{LN21},Theorem 4.6]\label{thm:LN}
Let $0\le \dd<1$ and suppose that the symbol $\sigma:\RR^n\times \NN_0^n\to \mathbb{C}$ satisfies $\sigma\in\SM^{0,\K,\N}_{1,\dd}$ for 
$$\K= n+1\qquad\text{and}\qquad  \N=2\Big\lceil\dfrac{n+1}{2(1-\dd)}\Big\rceil.$$
Then the pseudo-multiplier $\sigma(\cdot,\LL)$ extends to a bounded operator on $L^p(\RR^n)$ for every $1<p<\infty$. 
\end{thm}
\noindent For $\dd=0$ a similar set of conditions was given in \cite[Theorem 1.4]{BT}, but there the $L^2(\RR^n)$ boundedness was assumed. Theorem \ref{thm:LN} removes this assumption but requires more regularity in the spatial variable (only $\N=1$ is required in \cite[Theorem 1.4]{BT}).  

\bigskip

In view of \eqref{eq:CVpdo} an interesting question that arises is: {\it whether the number of derivatives in Theorem \ref{thm:LN} can be reduced?} The aim of the present article is to give additional answers to these questions.
Our main result is a new proof for case $\rr=1$ and $\dd<1$ in Theorem \ref{thm:LN} that requires less regularity in the frequency variable.

\begin{Theorem}\label{thm:main1}
Let $0\le \dd<1$ and suppose that the symbol $\sigma:\RR^n\times \NN_0^n\to \mathbb{C}$ satisfies $\sigma\in\SM^{0,\K,\N}_{1,\dd}$ for 
$$\K= \floor{n/2}+1\qquad\text{and}\qquad  \N= \Big\lfloor{\f{2n}{1-\dd}\Big\rfloor}+1.$$
Then the pseudo-multiplier $\sigma(\cdot,\LL)$ extends to a bounded operator  on $L^2(\RR^n)$. 
\end{Theorem}

\noindent On comparison with Theorem \ref{thm:LN} we see that the required number of derivatives in the frequency has been reduced from $n+1$ to $\floor{n/2}+1$, albeit at the cost of increased spatial regularity. 

We wish now to offer a few comments on our proof of Theorem \ref{thm:main1}. As in the original proof of the Calderon--Vaillancourt's result \eqref{eq:CVpdo}, our method relies on a version of the Cotlar--Knapp--Stein lemma (see Lemma \ref{lem:CKS}) but with some key differences. 
To highlight the differences, recall that the analysis of pseudo-differential operators \eqref{eq:pdo} on $\RR^n$ exploits an integration by parts in the frequency variables
$$|x-y|^2e^{i\xi\cdot x}\overline{e^{i\xi\cdot y}}=-\Delta_\xi e^{i\xi\cdot x}\overline{e^{i\xi\cdot y}},$$
and,  because of the symmetry in both variables,  an analogous integration by parts in space.
In the case of Hermite pseudo-multipliers however, due their asymmetry, one  requires  a different analysis in each  of the frequency and spatial variables. For summation in the frequency one can use a summation-by-parts type formula (see \eqref{eq:FDI}) that harks back to  \cite{T87,T93}. This type of formula produces an almost orthogonal decay on $\sum_\xi \sigma(x, \xi)h_\xi(x) h_\xi(y)$ in the spatial variables $|x-y|$, 
and has been recently exploited in \cite{BLL21,BT,LN21,LN22,PX} to further develop analysis for Hermite expansions. 

To fully exploit the Cotlar--Knapp--Stein lemma however, one needs a corresponding method to analyse the summation in the spatial variable in order to produce an almost orthogonal decay on $\ip{\sigma(\cdot,\xi)h_\xi,h_\eta}$ in terms of $|\xi-\eta|$. In certain cases explicit formulae for $\ip{\sigma(\cdot)h_\xi,h_\eta}$ are well known and classical. See for example \cite{Wang} for a treatment of the gaussian $e^{-x^2}$  on $\RR$, where one has essentially  an exponential decay in $|\xi-\eta|$ away from the diagonal and a polynomial decay  $\ip{e^{-x^2} h_\xi,h_\eta} \lesi (\xi-\eta)^{-1/2}$ near the diagonal. 

One of the contributions of the present article is the establishment of an integration-by-parts type formula that performs a similar role to the summation-by-parts in frequency \eqref{eq:FDI}, but for summation in the spatial variable (see Proposition \ref{prop:IBP}). The nucleus of this lemma is the classical Lagrange identity for orthogonal expansions, which can be elegantly used to show that $\{h_\xi\}_\xi$ are orthogonal. The insight behind Lemma  \ref{prop:IBP} is that through an extension of this argument one can derive almost orthogonality for expressions of the form $\ip{\sigma(\cdot,\xi)h_\xi,h_\eta}$. 

\bigskip

We give two additional results that require  regularity only in one of the variables in the symbol up to $\floor{n/2}+1$, both with relatively straightforward proofs. 

The first assumes regularity in the frequency but not in the spatial variable; such symbols are inspired by the `rough' symbols  $L^\infty S^m_\rr$ investigated in \cite{KS}.
\begin{Theorem}\label{thm:main2}
Let $m<0$ and suppose that $\sigma:\RR^n\times \NN_0^n\to \mathbb{C}$ satisfies 
\begin{align*}
\big\Vert\diff_\xi^\kappa \sigma(\cdot,\xi)\big\Vert_{L^\infty(\RR^n)}\le C_\kappa \ip{\xi}^{\f{m}{2}-|\kappa|},\qquad |\kappa|\le \floor{n/2}+1.
\end{align*}
Then the pseudo-multiplier $\sigma(\cdot,\LL)$ extends to a bounded operator  on $L^2(\RR^n)$. 
\end{Theorem}

As alluded to earlier, if $\sigma:\NN_0^n\to \mathbb{C}$ is a bounded function then  the multiplier $\sigma(\LL)$ is automatically bounded on $L^2(\RR^n)$. Our second result is a criterion along these lines and assumes regularity only in the spatial variable.
\begin{Theorem}\label{thm:main3}
Suppose the symbol $\sigma:\RR^n\times\NN_0^n\to\mathbb{C}$ satisfies 
\begin{align}\label{eq:sobolev-cond}
\sup_\xi\big\Vert\D_x^\nu \sigma(\cdot,\xi)\big\Vert_{L^2(\RR^n)}\le C_\nu,\qquad |\nu|\le \floor{n/2}+1.
\end{align}
Then the  Hermite pseudo-multiplier $\sigma(\cdot,\LL)$ extend to bounded operators  on  $L^2(\RR^n)$. 
\end{Theorem}
\noindent Note that this result is an analogue of \cite[Theorem 4.8.1]{RT} for periodic pseudo-differential operators. The main  insight here is that, unlike pseudo-differential operators on $\RR^n$, the orthogonality of Hermite expansions enables a similar argument to hold for Hermite pseudo-multipliers.

\begin{Remark}\label{rem:intro}
We conclude this section with several remarks on our results above. 
\begin{enumerate}[(i)]
\item It is worth noting  that  the proof of Theorem \ref{thm:main1} can be extended to $\rr<1$ and $m=n(\rr-1)$. In this case the required number of derivatives are $\K=\floor{n/2}+1$ and $\N=\floor{2n\rr/(1-\dd)}+1$.  However for brevity we only focus here on the case $\rho=1$  and leave the details for $\rr<1$ to the interested reader. 
\item Theorem \ref{thm:main3} also holds when the derivative $\partial$ is replaced by other kinds of derivatives such as $A^*$ (see Section \ref{sec:prelim}). In this case, one needs a Sobolev-type embedding for $A^*$, which can be found in \cite[Theorem 4 and Theorem 8(ii)]{BTo}. For other similar $L^2$ Sobolev-type conditions on the symbol we refer the reader to \cite[Corollary 5.6]{LN21} (for Hermite pseudo-multipliers) and \cite[Corollary 2.2]{Hwang} (for pseudodifferential operators on $\RR^n$).
\item Regarding results for $L^p(\RR^n)$ for $p\ne 2$, we note that the operators $\sigma(\cdot,\LL)$ in Theorem \ref{thm:LN} are actually Calder\'on--Zygmund operators. For $\dd=0$ this was essentially proved  in  \cite[page 165]{BT}. We also note that Theorem \ref{thm:main2} can also extended to $p\ne 2$. These facts will appear in forthcoming work.
\end{enumerate}
\end{Remark}

\subsection{Pseudo-multipliers associated with Hermite polynomials}

The second type of expansion we wish to consider are those defined using \emph{normalized Hermite polynomial} expansions $\HH_\xi$, which are defined by 
\begin{align}\label{eq:hermitepoly}
\HH_\xi(x) = \big( 2^{|\xi|} \xi !\big)^{-\half} H_\xi(x),
\end{align}
where $H_\xi(x):=\prod_{j=1}^nH_{\xi_j}(x_j)$ are the $n$-dimensional Hermite polynomials. 
The collection $\{\HH_\xi\}_{\xi\in\NN_0^n}$ forms a complete orthonormal basis for $L^2(\gamma)=L^2(\RR^n,d\gamma)$ with respect to the inner product given by $\ip{f,g}_\gamma :=\int_{\RR^n} fg\,d\gamma$, where $d\gamma=\pi^{-\f{n}{2}}e^{-|x|^2}dx$ is the Gaussian measure. See \cite{Ur} for further details. 

Given a symbol $\sigma:\RR^n\times \NN_0^n\to \mathbb{C}$ we define the `Gaussian pseudo-multiplier'  by
\begin{align}\label{eq:gpdo}
\sigma(x,\OU) f(x):= \sum_{\xi\in\NN_0^n} \sigma(x,\xi) \ip{f,\HH_\xi}_\gamma \HH_\xi(x),
\end{align}
where $\OU:=-\half \Delta +x\cdot\nabla$ is the Ornstein--Uhlenbeck operator. 

As in the Hermite function case,  when $\sigma$ is independent of $x$ (that is, $\sigma(\cdot,\OU)=\sigma(\OU)$ is a spectral multiplier), then through a similar argument one can see that $\sigma\in L^\infty(\RR^n)$ is a necessary and sufficient for the $L^2(\gamma)$-boundedness of $\sigma(\OU)$ (see \cite[Chapter 6]{Ur}).

While spectral multipliers and other operators in this setting have been thoroughly explored (see \cite{Ur} for an extensive survey), the study of pseudo-differential type operators is yet in its infancy; indeed it was remarked in \cite{Ur} that operators of the form \eqref{eq:gpdo} `require further analysis and study'. We contribute to this line of inquiry by giving sufficient conditions for the $L^2$-boundedness of pseudo-multipliers \eqref{eq:gpdo}.

In fact, it turns out that there is a close link between pseudo-multipliers in the Gaussian \eqref{eq:gpdo} and Hermite \eqref{eq:hpdo} settings, allowing us to transfer results from one to the other. In \cite{AT} this transferability is expressed by the statement that `operators associated to $\LL$ and $\OU$ are unitarily equivalent in $L^2$'.

More precisely, for pseudo-multipliers we have the following:
given a symbol $\sigma :\RR^n\times \NN_0^n\to \mathbb{C}$, then $\sigma(\cdot,\OU)$ is bounded on $L^2(\gamma)$ if and only if $\sigma(\cdot,\LL)$ is bounded on $L^2(\RR^n)$; furthermore,
\begin{align}\label{eq:hermite-equiv}
\Vert \sigma(\cdot,\OU) \Vert_{L^2(\gamma)\to L^2(\gamma)} = \Vert \sigma(\cdot,\LL)\Vert_{L^2\to L^2}.
\end{align}
To see \eqref{eq:hermite-equiv}, consider the operator $\U$ defined by $\U f(x):=\pi^{-n/4} e^{-|x|^2/2}f(x)$ for any function $f$. Then $\U$ is an isometry from $L^2(\RR^n)$ into $L^2(\gamma)$; that is,  we have $\Vert \U f\Vert_{L^2} = \Vert f\Vert_{L^2(\gamma)}$ (\cite[Lemma 3.1]{AT}). In addition, it can be readily seen from the definitions that $\U \HH_\xi = h_\xi$ and  $\ip{f,\HH_\xi}_\gamma = \ip{\U f, h_\xi}$. These facts, along with  some simple calculations give
$$ \Vert \sigma(\cdot,\OU)\Vert_{L^2(\gamma)\to L^2(\gamma)} = \big\Vert \U^{-1}\circ\sigma(\cdot,\LL)\circ\U\big\Vert_{L^2(\gamma)\to L^2(\gamma)}.$$
Since $\U$ is an isometry, then \eqref{eq:hermite-equiv} follows. 

With \eqref{eq:hermite-equiv} in hand, we immediately have the following.
\begin{Corollary}\label{cor:gpdo}
Suppose that the symbol $\sigma:\RR^n\times\NN_0^n\to\mathbb{C}$ satisfies any  one of the conditions in Theorem \ref{thm:LN}, Theorem \ref{thm:main1}, Theorem \ref{thm:main2} or Theorem \ref{thm:main3}. Then the Gaussian pseudo-multiplier given by \eqref{eq:gpdo} extends to a bounded operator on $L^2(\gamma)$. 
\end{Corollary}

\noindent Corollary \ref{cor:gpdo} furnishes some sufficient conditions of $L^2$ boundedness of Gaussian pseudo-multipliers. The question of $L^p$ boundedness remains open and we shall reserve this for future studies.

\subsection{Further comments}

While this work was in its final stages, the author came across the interesting work of \cite{BG} which gives a Calder\'on--Vaillancourt theorem for Hermite pseudo-multipliers. In \cite[Theorem 1.3]{BG} the authors show that  if $\sigma\in\SM^{0,\infty,\infty}_{\rr,\dd}$ for $0\le \dd/2\le \rr\le 1$, $\dd\ne 2$,  then the Hermite pseudo-multiplier $\sigma(\cdot,\LL)$ is bounded on $L^2(\RR^n)$. 
In comparison with  Theorem \ref{thm:main1} we see that both results overlap when $\rr=1$, but otherwise deal with different regions (see Remark \ref{rem:intro} (i)). In addition,  the techniques employed here are different and the proof of \cite[Theorem 1.3]{BG} is not given explicitly but is rather described as a consequence of a similar result for the Grushin operator. In the present article, we work directly in the Hermite setting which we hope will be of value to future researchers in the topic.

Finally, it is worth pointing out that  \cite[Theorem 1.3]{BG} along with \eqref{eq:hermite-equiv} establishes a Calder\'on--Vaillancourt result for Gaussian pseudo-multipliers.

\medskip
{\bf Organization:}
The organization of this article is as follows. 
In Section \ref{sec:prelim} we give some additional necessary notation and background on Hermite functions and the associated operators. We also describe and prove some important ingredients and tools needed in the proofs of our main results, including summation and integration by parts in the Hermite context (Section \ref{sec:byparts}) and a Littlewood--Paley type decomposition and associated kernel estimates (Section \ref{sec:decomp}). With these preliminaries in hand, we give the proofs of Theorems \ref{thm:main1} - \ref{thm:main3} in Section \ref{sec:L2}.

\medskip
{\bf Acknowledgements:}  The author thanks Virginia Naibo and Rahul Garg for  helpful discussions and comments. The author also extends their gratitude to the referee for their advice and suggestions which has improved the quality of the paper.

\section{Preliminaries}\label{sec:prelim}
In this section  we describe  some notation and basic facts related to Hermite functions. We also give some important tools and estimates that will be used in the proofs in rest of the paper. These include integration and summation by-parts in the context of Hermite expansions (Section \ref{sec:byparts}) and a decomposition of symbols and their relevant kernel estimates (Section \ref{sec:decomp}).

For each $i\in\{1,\dots,n\}$ the Hermite derivatives are defined by
\begin{align*}
A_i^* = \f{\D}{\D x_i}+x_i \qquad\text{and}\qquad A_i=-\f{\D}{\D x_i}+x_i.
\end{align*}
These are sometimes called the \emph{annihilation} and \emph{creation} operators respectively. Moreover the operator $\LL$ can be factored as $\LL= \f{1}{2} \sum_i^n(A_iA_i^* + A_i^*A_i)$. Additional properties relevant to this paper is described in Lemmas \ref{lem:FDItools} below. For further details and background concerning the Hermite functions and associated operators the reader may consult \cite{T93}.

For each $N\in\NN_0$  we  define the orthogonal projection of $f$ onto $\bigoplus^N_{k=0} \text{span}\{h_\xi: |\xi|=k\}$ by
\begin{align}\label{eq:QQ}
	\QQ_N f = \sum_{k=0}^N \sum_{|\xi|=k} \ip{f, h_\xi} h_\xi
	\qquad\text{with kernel}\qquad
	\QQ_N(x,y)=\sum_{k=0}^N \sum_{|\xi|=k}h_\xi(x)h_\xi(y).
\end{align}
The following bounds are known (see \cite[p.376]{PX}): there exists $\vartheta>0$ such that for any $N\in \NN$
\begin{align}\label{QQ bound}
	\QQ_{N}(x,x) 
	\lesssim \left
\lbrace 
	\begin{array}{ll}
			N^{n/2}  \qquad &\forall x,\\
			e^{-2\vartheta|x|^2} \qquad &\text{if}\quad |x| \ge \sqrt{4N+2}.
				\end{array}
\right.
\end{align}

\subsection{Summation and integration by-parts for Hermite function expansions}\label{sec:byparts}
In this section we list two key identities for hermite expansions that are crucial for the proofs of our main results. 

The first can be viewed as an integration-by-parts in the frequency variable and has its origins in \cite{T87,T93}. It is used to obtain the kernel estimates of Lemma \ref{lem:Kj}.
\begin{Proposition}[By-parts in frequency variables]\label{prop:FDI}
Suppose $k$ is a function defined on $\RR^n\times\RR^n\times\NN_0^n$ and set
$$ K(x,y) = \sum_{\xi\in\NN_0^n} k(x,y,\xi) \,h_\xi(x)h_\xi(y).$$
If $N\in\ZZ_+,$ it holds that
\begin{align}\label{eq:FDI}
	2^N(x_i- y_i)^N K(x,y) = \sum_{\substack{\xi\in\NN_0^n\\N/2\le \ell\le N,\\ \nu,\om\ge 0\\ \nu+\om=2\ell-N}} c_{\nu,\ell, N} \,d_\nu(\xi_i) d_\om(\xi_i)\, \diff_i^\ell k(x,y,\xi)h_{\xi+\nu e_i}(x) h_{\xi+\om e_i}(y),
\end{align}
where $c_{\nu,\ell,N} = (-1)^{\ell-\nu}4^{N-\ell}(2N-2\ell-1)!!\binom{N}{2\ell-N}$ and for $\lambda \ge0$, $d_m(\lambda) = \prod_{r=0}^{m-1}\sqrt{2(\lambda+r)+2}$ if $m\ge 1$ and $d_m(\lambda) =1$ if $m=0$.
\end{Proposition}
\begin{Remark}\label{rem:FDI}
Observe that the quantities $d_m(\lambda)$ satisfies $ d_{m}(\lambda)\lesi \ip{\lambda}^{m/2}$. This fact will be  frequently employed in the sequel. 
\end{Remark}
\noindent One can trace \eqref{eq:FDI} back to \cite{T93}; the integer version appears in \cite[Lemma 3.2.3]{T93} and a multi-index version appears in  \cite[(4.2.12)]{T93}. The constants were not calculated explicitly there, but this was done in  \cite[Lemma 8]{PX} for the integer version. The identity in Proposition \ref{prop:FDI} is the multi-index version of \cite[Lemma 8]{PX} and we give the proof here for completeness.

\bigskip
Our second result is an integration-by-parts in the spatial variable and is new. It is used in the proof of  Theorem \ref{thm:main1}.

\begin{Proposition}[By-parts in spatial variables]\label{prop:IBP}
Let $N\in\NN_0$ and  $g:\RR^n \times \NN_0^n\times\NN_0^n\to\mathbb{C}$ be a $C^N$ function in the first variable such that $|\D_x^\gamma g(x,\xi,\eta)h_\xi(x)|\to 0$ as $|x|\to \infty$ for all $|\gamma|\le N$ and $(\xi,\eta)\in\NN_0^{2n}$.
Set $$G(\xi,\eta) := \int g(x,\xi,\eta) \,h_\xi(x)\,h_\eta(x)\,dx $$
Then for each $i\in\{1,\dots,n\}$, 
\begin{align}\label{eq:IBP}
&2^N(\xi_i-\eta_i)^N G(\xi,\eta) \\
&\qquad\qquad=\sum_{\substack{-N\le \alpha,\beta\le N,\\ 1\le \ell\le N,\\ |\alpha-\beta|\le N}} C(N,\xi_i,\eta_i,\alpha,\beta)\int \D_i^\ell g(x,\xi,\eta) \,h_{\xi+\alpha e_i}(x) \,h_{\eta+\beta e_i}(x)\,dx 
\notag
\end{align}
where $C(N,\xi_i,\eta_i,\alpha,\beta)$ are constants satisfying
\begin{align}\label{eq:IBP2}
|C(N,\xi_i,\eta_i,\alpha,\beta)| \le c_N \ip{\xi_i \vee \eta_i}^{N/2},
\end{align}
and the sum contains at most $5^N$ terms. 
\end{Proposition}

\subsubsection{Proof of Proposition \ref{prop:FDI}}
The result will be a consequence of the following two lemmas.

\begin{Lemma}\label{lem:FDItools}
\begin{align}
A_i h_\xi &= \sqrt{2\xi_i +2}\,h_{\xi+e_i},\qquad \text{where} \quad A_i =-\partial_i+x_i.\label{FDI1}\\
A^*_i h_\xi &= \sqrt{2\xi_i}\,h_{\xi-e_i},\qquad \text{where} \quad A^*_i =\partial_i+x_i.\label{FDI2}
\end{align}
\begin{align}\label{FDI3}
2x_i h_\xi =\sqrt{2\xi_i+2}\, h_{\xi+e_i} +\sqrt{2\xi_i}\,h_{\xi-e_i}.
\end{align}
\begin{align}\label{FDI4}
(x_i-y_i)\big(\A{y}_i-\A{x}_i\big)^r=\big(\A{y}_i-\A{x}_i\big)^r(x_i-y_i)-2r\big(\A{y}_i-\A{x}_i\big)^{r-1}, \qquad r\ge1.
\end{align}
\begin{align}\label{FDI5}
\big(\A{y}_i-\A{x}_i\big)\,h_\xi(x)\,h_\xi(y)=\sqrt{2\xi_i+2}\,\big(h_\xi(x)h_{\xi+e_i}(y)-h_{\xi+e_i}(x)h_\xi(y)\big).
\end{align}
\begin{align}\label{FDI6}
2(x_i-y_i)h_\xi(x)h_\xi(y)=\left\lbrace \begin{array}{cl}
    - \big(\A{y}_i-\A{x}_i\big)\diff_i [h_{\xi-e_i}(x) h_{\xi-e_i}(y)] &\quad \text{if}\quad \xi_i\ge 1,\\
    -\big(\A{y}_i-\A{x}_i\big)h_\xi(x)h_\xi(y) &\quad \text{if}\quad\xi_i=0.
 \end{array}\right. 
\end{align}
For any $k\ge 0$, 
\begin{align}\label{FDI7}
&\sum_{\xi\in\NN_0^n}f(x,y,\xi)\big(\A{y}_i-\A{x}_i\big)^k2(x_i-y_i)h_\xi(x)h_\xi(y)\\
&\hspace{100pt}=\sum_{\xi\in\NN_0^n}\diff_if(x,y,\xi)\big(\A{y}_i-\A{x}_i\big)^{k+1}h_\xi(x)h_\xi(y). \notag
\end{align}
\end{Lemma}

\begin{Lemma}\label{lem:FDIa}
For each $N\in\ZZ_+$ we have
\begin{align}\label{eq:FDIa}
	2^N(x_i- y_i)^N \KK(x,y) = \sum_{\f{N}{2}\le \ell\le N}c_{\ell, N}\sum_{\xi\in \NN_0^n} \diff_i^\ell k(x,y,\xi)\big(\A{y}_i-\A{x}_i\big)^{2\ell-N} \,h_\xi(x)\,h_\xi(y)
\end{align}
where $c_{\ell,N} = (-4)^{N-\ell}(2N-2\ell-1)!!\binom{N}{2\ell-N}$.
\end{Lemma}

Let us give  the proof of Proposition \ref{prop:FDI} assuming the above two results. Firstly by applying \eqref{FDI1} repeatedly we have easily
$$ (\A{y}_i)^\nu h_\xi(y) = \prod_{r=0}^{\nu-1} \sqrt{2(\xi_i+r)+2} h_{\xi+\nu e_i}(y)=d_\nu(\xi_i)h_{\xi+\nu e_i}(y).$$
Then using the binomial theorem we have
\begin{align*}
\big(\A{y}_i-\A{x}_i\big)^{2\ell-M}h_\xi(x)h_\xi(y) 
&=\sum_{\substack{\nu,\om\ge0\\ \nu+\om=2\ell-N}}  (-1)^\om \tbinom{2\ell-M}{\nu,\om}(\A{y}_i)^\nu h_\xi(y) (\A{x}_i)^{\om} h_\xi(x) \\
&= \sum_{\substack{\nu,\om\ge0\\ \nu+\om=2\ell-N}}  (-1)^\om \tbinom{2\ell-M}{\nu,\om} d_\nu(\xi_i)d_\om(\xi_i) h_{\xi+\nu e_i}(y)  h_{\xi+\om e_i}(x) 
\end{align*}
Inserting this equation into \eqref{eq:FDIa} gives \eqref{eq:FDI}.

It remains now to give the proofs of Lemmas \ref{lem:FDItools} and \ref{lem:FDIa}. 
\begin{proof}[Proof of Lemma \ref{lem:FDItools}]
The first two identities are standard. Combining  \eqref{FDI1} and \eqref{FDI2} we get \eqref{FDI3}. The identity \eqref{FDI4} can be found in \cite[(3.2.23)]{T93}, while identity \eqref{FDI5} follows from repeated application of \eqref{FDI1}. Let us prove \eqref{FDI6}. For $\xi_i\ge 1$ we have by \eqref{FDI5} and \eqref{FDI3}
\begin{align*}
& - \big(\A{y}_i-\A{x}_i\big)\diff_i [h_{\xi-e_i}(x) h_{\xi-e_i}(y)] \\
 &\qquad= \big(\A{y}_i-\A{x}_i\big)[h_{\xi-e_i}(x) h_{\xi-e_i}(y) -h_\xi(x)h_\xi(y)]  \\
 &\qquad=\sqrt{2\xi_i}\,\big(h_{\xi-e_i}(x)h_{\xi}(y)-h_{\xi}(x)h_{\xi-e_i}(y)\big) 
 -\sqrt{2\xi_i+2}\,\big(h_{\xi}(x)h_{\xi+e_i}(y)-h_{\xi+e_i}(x)h_{\xi}(y)\big)  \\
 &\qquad=h_\xi(y)\big(\sqrt{2\xi_i}h_{\xi-e_i}(x)+\sqrt{2\xi_i+2}h_{\xi+e_i}(x)\big)-h_\xi(x)\big(\sqrt{2\xi_i}h_{\xi-e_i}(y)+\sqrt{2\xi_i+2}h_{\xi+e_i}(y)\big)\\
  &\qquad=h_\xi(y)\big(2x_ih_{\xi}(x)\big)-h_\xi(x)\big(2y_ih_{\xi}(y)\big) \\
  &\qquad= 2(x_i-y_i)h_\xi(x)h_\xi(y).
\end{align*}
Using the fact that in one dimensions we have
$$ h_1(t)=\f{2}{\sqrt{2}}t h_0(t),$$
then we obtain
\begin{align*}
h_{\xi+e_i}(x)
=\Big[\prod^n_{j\ne i}h_{\xi_j}(x_j)\Big] \times \f{2}{\sqrt{2}}x_ih_0(x_i)
=\f{2}{\sqrt{2}}x_i \Big[\prod^n_{j=1}h_{\xi_j}(x_j)\Big] 
=\f{2}{\sqrt{2}}x_i h_\xi(x).
\end{align*}
This fact along with  \eqref{FDI5}  gives, for $\xi_i=0$,
\begin{align*}
 -\big(\A{y}_i-\A{x}_i\big)h_\xi(x)h_\xi(y)
 &=-\sqrt{2}\big(h_\xi(x)h_{\xi+e_i}(y)-h_{\xi+e_i}(x)h_\xi(y)\big) \\
 &=-\sqrt{2}\big(h_\xi(x)\big(\f{2}{\sqrt{2}}y_ih_\xi(y)\big) -h_\xi(y)\big(\f{2}{\sqrt{2}}x_ih_\xi(x)\big)\big)\\
 &=2(x_i-y_i)h_\xi(x)h_\xi(y).
\end{align*}
We now turn to \eqref{FDI7}. Firstly, reindexing  $\xi$ by $\xi-e_i$ gives
\begin{align*}
&\sum_{\xi\in\NN_0^n}\diff_if(x,y,\xi)\big(\A{y}_i-\A{x}_i\big)^{k+1}h_\xi(x)h_\xi(y)\\
&\hspace{80pt}=\sum_{\xi\in\NN_0^n}f(x,y,\xi+e_i)\big(\A{y}_i-\A{x}_i\big)^{k+1}h_\xi(x)h_\xi(y)\\
&\hspace{140pt}-\sum_{\xi\in\NN_0^n}f(x,y,\xi)\big(\A{y}_i-\A{x}_i\big)^{k+1}h_\xi(x)h_\xi(y)\\
&\hspace{80pt}=\sum_{\xi_i\ge1}f(x,y,\xi)\big(\A{y}_i-\A{x}_i\big)^{k+1}h_{\xi-e_i}(x)h_{\xi-e_i}(y)\\
&\hspace{140pt}-\sum_{\xi\in\NN_0^n}f(x,y,\xi)\big(\A{y}_i-\A{x}_i\big)^{k+1}h_\xi(x)h_\xi(y). 
\end{align*}
Now separating $\xi_i=0$ out of the second sum, we obtain
\begin{align*}
&\sum_{\xi\in\NN_0^n}\diff_if(x,y,\xi)\big(\A{y}_i-\A{x}_i\big)^{k+1}h_\xi(x)h_\xi(y)\\
&\hspace{80pt}=\sum_{\xi_i\ge1}f(x,y,\xi)\big(\A{y}_i-\A{x}_i\big)^{k+1}h_{\xi-e_i}(x)h_{\xi-e_i}(y)\\
&\hspace{140pt}-\sum_{\xi_i\ge 1}f(x,y,\xi)\big(\A{y}_i-\A{x}_i\big)^{k+1}h_\xi(x)h_\xi(y)\\
&\hspace{180pt}-\sum_{\xi_i=0}f(x,y,\xi)\big(\A{y}_i-\A{x}_i\big)^{k+1}h_\xi(x)h_\xi(y)  \\
&\hspace{80pt}=\sum_{\xi_i\ge1}f(x,y,\xi)\big(\A{y}_i-\A{x}_i\big)^{k+1}(-\diff_i)\big[ h_{\xi-e_i}(x)h_{\xi-e_i}(y)\big]\\
&\hspace{140pt}-\sum_{\xi_i=0}f(x,y,\xi)\big(\A{y}_i-\A{x}_i\big)^{k+1}h_\xi(x)h_\xi(y).
\end{align*}
Finally, applying  \eqref{FDI6} we get
\begin{align*}
&\sum_{\xi\in\NN_0^n}\diff_if(x,y,\xi)\big(\A{y}_i-\A{x}_i\big)^{k+1}h_\xi(x)h_\xi(y)\\
&\hspace{80pt}=\sum_{\xi_i\ge 1}f(x,y,\xi)\big(\A{y}_i-\A{x}_i\big)^k2(x_i-y_i)h_\xi(x)h_\xi(y)\\
&\hspace{140pt}+\sum_{\xi_i=0}f(x,y,\xi)\big(\A{y}_i-\A{x}_i\big)^k2(x_i-y_i)h_\xi(x)h_\xi(y)\\
&\hspace{80pt}=\sum_{\xi\in\NN_0^n}f(x,y,\xi)\big(\A{y}_i-\A{x}_i\big)^k2(x_i-y_i)h_\xi(x)h_\xi(y),
\end{align*}
which gives \eqref{FDI7}.
\end{proof}

\begin{proof}[Proof of Lemma \ref{lem:FDIa}]
We shall prove the Lemma by induction on $N$. The case $N=1$,
\begin{align*}
2(x_i-y_i)\KK(x,y)=\sum_{\xi\in\NN_0^n}\diff_i k(x,y,\xi)(\A{y}_i-\A{x}_i)h_\xi(x) h_\xi(y),
\end{align*}
follows immediately by applying \eqref{FDI7} with $f(x,y,\xi)=k(x,y,\xi)$ and $k=0$.

Let us turn to the inductive step; note that $\ell\ge\f{N+1}{2}$. We write
\begin{align*}
&2^{N+1}(x_i- y_i)^{N+1} \KK(x,y) \\
&\qquad=2(x_i-y_i) \big[  2^{N}(x_i- y_i)^{N} \KK(x,y) \big]\\
&\qquad=2(x_i-y_i)\sum_{\f{N}{2}\le \ell\le N}c_{\ell, N}\sum_{\xi\in \NN_0^n} \diff_i^\ell k(x,y,\xi)\big(\A{y}_i-\A{x}_i\big)^{2\ell-N} \,h_\xi(x)\,h_\xi(y) \\
&\qquad=\sum_{\f{N}{2}\le \ell\le N}c_{\ell, N}\sum_{\xi\in \NN_0^n} \diff_i^\ell k(x,y,\xi)\Big[2(x_i-y_i)\big(\A{y}_i-\A{x}_i\big)^{2\ell-N} \,h_\xi(x)\,h_\xi(y)\Big] \\
&\qquad=\sum_{\f{N}{2}\le \ell\le N}c_{\ell, N}\sum_{\xi\in \NN_0^n} \diff_i^\ell k(x,y,\xi)\big(\A{y}_i-\A{x}_i\big)^{2\ell-N} 2(x_i-y_i)\,h_\xi(x)\,h_\xi(y) \\
&\qquad\qquad-\sum_{\f{N}{2}\le \ell\le N}c_{\ell, N}\sum_{\xi\in \NN_0^n} \diff_i^\ell k(x,y,\xi)4(2\ell-N)\big(\A{y}_i-\A{x}_i\big)^{2\ell-N-1} \,h_\xi(x)\,h_\xi(y) \\
&\qquad =:I+II
\end{align*}
by \eqref{FDI4} with $r=2\ell-N \ge 1$ in the last step. 

Now applying \eqref{FDI7} with $f(x,y,\xi)=\diff_i^\ell k(x,y,\xi)$ and $k=2\ell-N$, and then reindexing $\ell$ to $\ell-1$, we have
\begin{align*}
I &=\sum_{\f{N}{2}\le \ell\le N}c_{\ell, N}\sum_{\xi\in \NN_0^n} \diff_i^{\ell+1} k(x,y,\xi)\big(\A{y}_i-\A{x}_i\big)^{2\ell-N+1} \,h_\xi(x)\,h_\xi(y) \\
&=\sum_{\f{N}{2}\le \ell-1\le N}c_{\ell-1, N}\sum_{\xi\in \NN_0^n} \diff_i^\ell k(x,y,\xi)\big(\A{y}_i-\A{x}_i\big)^{2\ell-N-1} \,h_\xi(x)\,h_\xi(y)\\
&=c_{N, N}\sum_{\xi\in \NN_0^n} \diff_i^{N+1} k(x,y,\xi)\big(\A{y}_i-\A{x}_i\big)^{N+1} \,h_\xi(x)\,h_\xi(y) \\
&\qquad + \sum_{\f{N+1}{2}\le \ell\le N}c_{\ell-1, N}\sum_{\xi\in \NN_0^n} \diff_i^\ell k(x,y,\xi)\big(\A{y}_i-\A{x}_i\big)^{2\ell-N-1} \,h_\xi(x)\,h_\xi(y).
\end{align*}
Here $c_{\ell,N}:=0$ whenever $\ell <N/2$, and note that the second sum is empty whenever $(N+1)/2\le \ell <(N+2)/2$. Next we have
\begin{align*}
II
&=\sum_{\f{N}{2}\le \ell\le N}c_{\ell, N}\sum_{\xi\in \NN_0^n} \diff_i^\ell k(x,y,\xi)(-4(2\ell-N))\big(\A{y}_i-\A{x}_i\big)^{2\ell-N-1} \,h_\xi(x)\,h_\xi(y)  \\
&= \sum_{\f{N+1}{2}\le \ell\le N} \big[-4(2\ell-N)c_{\ell, N}\big]\sum_{\xi\in \NN_0^n} \diff_i^\ell k(x,y,\xi)\big(\A{y}_i-\A{x}_i\big)^{2\ell-N-1} \,h_\xi(x)\,h_\xi(y),
\end{align*}
again noting that the sum is empty for $N/2\le \ell < (N+1)/2$.

Combining the above terms $I$ and $II$ we obtain
\begin{align*}
&2^{N+1}(x_i- y_i)^{N+1} \KK(x,y) \\
&\hspace{20pt}=\sum_{\f{N+1}{2}\le \ell\le N}\big[c_{\ell-1, N}-4(2\ell-N)c_{\ell, N}\big]\sum_{\xi\in \NN_0^n} \diff_i^\ell k(x,y,\xi)\big(\A{y}_i-\A{x}_i\big)^{2\ell-N-1} \,h_\xi(x)\,h_\xi(y) \\
&\hspace{60pt}+c_{N, N}\sum_{\xi\in \NN_0^n} \diff_i^{N+1} k(x,y,\xi)\big(\A{y}_i-\A{x}_i\big)^{N+1} \,h_\xi(x)\,h_\xi(y) \\
&\hspace{20pt}= \sum_{\f{N+1}{2}\le \ell\le N+1}c_{\ell,N+1}\sum_{\xi\in \NN_0^n} \diff_i^\ell k(x,y,\xi)\big(\A{y}_i-\A{x}_i\big)^{2\ell-N-1} \,h_\xi(x)\,h_\xi(y)
\end{align*}
where 
$$ c_{\ell,N+1}=\left\lbrace \begin{array}{cl}
    c_{N,N} &\quad \text{if}\quad \ell=N+1,\\
    c_{\ell-1, N}-4(2\ell-N)c_{\ell, N} &\quad \text{if} \quad\f{N+1}{2}\le \ell\le N.
 \end{array}\right. 
 $$
At this point the constants can be obtained by solving the recursion as in \cite[Lemma 8]{PX}.
\end{proof}

\subsubsection{Proof of Proposition \ref{prop:IBP}}
We first derive the following identity which will be crucial for our proof. Given $g$ as defined in Proposition \ref{prop:IBP}, we have, for all $\xi,\eta\in\NN_0^n$ and $\alpha,\beta\in \ZZ$ with $\xi_i+\alpha\ge0$ and $\eta_i+\beta\ge0$,
\begin{align}\label{eq:lagrange2}
&2(\xi_i-\eta_i)\int g(x,\xi,\eta) h_{\xi+\alpha e_i}(x) h_{\eta+\beta e_i}(x)\,dx  \\
&\hspace{60pt}= \int \D_i g(x,\xi,\eta)\big[ \D_i h_{\xi+\alpha e_i}(x) h_{\eta+\beta e_i}(x)-h_{\xi+\alpha e_i}(x)\D_i h_{\eta+\beta e_i}(x)\big]dx  \notag\\
&\hspace{120pt} + 2(\beta-\alpha) \int g(x,\xi,\eta)h_{\xi+\alpha e_i}(x) h_{\eta+\beta e_i}(x)dx.  \notag
\end{align}
Equation \eqref{eq:lagrange2} is a consequence of 
Lagrange's classical identity  for Hermite expansions: 
\begin{align}\label{eq:lagrange}
v\LL_i u -u\LL_i v = \D_i (u\D_i v-v\D_i u),
\end{align}
where $\LL_i = -\D_i^2+x_i^2$.
Indeed, using the fact that $\LL_i h_\xi = (2\xi_i+1)h_\xi$, identity \eqref{eq:lagrange} and integration by parts we have
\begin{align*}
&2(\xi_i-\eta_i)\int g(x,\xi,\eta) h_{\xi+\alpha e_i}(x) h_{\eta+\beta e_i}(x)\,dx  \\
&\hspace{50pt}= \int g(x,\xi,\eta) \big[\LL_i h_{\xi+\alpha e_i}(x) h_{\eta+\beta e_i}(x)\,dx 
- h_{\xi+\alpha e_i}(x) \LL_i h_{\eta+\beta e_i}(x)\big]\,dx \\
&\hspace{100pt} +2(\beta-\alpha)\int g(x,\xi,\eta) h_{\xi+\alpha e_i}(x) h_{\eta+\beta e_i}(x)\,dx\\
&\hspace{50pt}= \int g(x,\xi,\eta) \D_i\big[ h_{\xi+\alpha e_i}(x) \D_i h_{\eta+\beta e_i}(x)\,dx 
- \D_i h_{\xi+\alpha e_i}(x) h_{\eta+\beta e_i}(x)\big]\,dx \\
&\hspace{100pt} +2(\beta-\alpha)\int g(x,\xi,\eta) h_{\xi+\alpha e_i}(x) h_{\eta+\beta e_i}(x)\,dx \\
&\hspace{50pt}= -\int \D_i g(x,\xi,\eta) \big[ h_{\xi+\alpha e_i}(x) \D_i h_{\eta+\beta e_i}(x)\,dx 
- \D_ih_{\xi+\alpha e_i}(x)  h_{\eta+\beta e_i}(x)\big]\,dx \\
&\hspace{100pt} +2(\beta-\alpha)\int g(x,\xi,\eta) h_{\xi+\alpha e_i}(x) h_{\eta+\beta e_i}(x)\,dx,
\end{align*}
which proves \eqref{eq:lagrange2}.

We can now continue with the proof of Proposition \ref{prop:IBP}. We will proceed by induction on $N$. In the sequel we set $a_k:=\sqrt{k/2}$ for $k\in\NN_0$. 
For $N=1$ we apply \eqref{eq:lagrange2} with $\alpha=\beta=0$. Then 
\begin{align*}
&2(\xi_i-\eta_i)\int g(x,\xi,\eta) \,h_\xi(x)\,h_\eta(x)\,dx \\
&\qquad = \int \D_i g(x,\xi,\eta) \big[\D_i h_\xi(x)h_\eta(x)-h_\xi(x)\D_i h_\eta(x)\big]\,dx \\
&\qquad = -a_{\xi_i+1}\int \D_i g(x,\xi,\eta) h_{\xi+e_i}(x)h_\eta(x)dx + a_{\xi_i}\int \D_i g(x,\xi,\eta) h_{\xi-e_i}(x)h_\eta(x)dx \\
&\qquad\qquad +a_{\eta_i+1}\int\D_i g(x,\xi,\eta) h_\xi(x) h_{\eta+e_i}(x)dx -a_{\eta_i}\int \D_ig(x,\xi,\eta)h_\xi(x)h_{\eta-e_i}(x)dx.
\end{align*}
Thus \eqref{eq:IBP} holds with $|C(1,\xi_i,\eta_i,0,0)|\lesi (\eta_i\vee \xi_i)^{1/2}$.

We now show the inductive step. Firstly, applying \eqref{eq:lagrange2} to the function $\D_i^\ell g$ we have, for any $\ell\in\NN_0$,
\begin{align*}
&2(\xi_i-\eta_i)\int \D_i^\ell g(x,\xi,\eta) h_{\xi+\alpha e_i}(x) h_{\eta+\beta e_i}(x)\,dx  \\
&\hspace{50pt}= \int \D_i^{\ell+1} g(x,\xi,\eta)\big[ \D_i h_{\xi+\alpha e_i}(x) h_{\eta+\beta e_i}(x)-h_{\xi+\alpha e_i}(x)\D_i h_{\eta+\beta e_i}(x)\big]dx \\
&\hspace{100pt} + 2(\beta-\alpha) \int \D_i^\ell g(x,\xi,\eta)h_{\xi+\alpha e_i}(x) h_{\eta+\beta e_i}(x)dx.
\end{align*}
Then, using this identity we obtain, for $N\ge1$,
\begin{align*}
&2^{N+1}(\xi_i-\eta_i)^{N+1}\int g(x,\xi,\eta) \,h_\xi(x)\,h_\eta(x)\,dx\\
&\hspace{40pt} =2(\xi_i-\eta_i)\sum_{\alpha,\beta,\ell,N} C(N,\xi_i,\eta_i,\alpha,\beta) \int \D_i^\ell g(x,\xi,\eta) \,h_{\xi+\alpha e_i}(x) \,h_{\eta+\beta e_i}(x)\,dx  \\
&\hspace{40pt}=\sum_{\alpha,\beta,\ell,N} C(N,\xi_i,\eta_i,\alpha,\beta) \\
&\hspace{80pt}\times\Big\{ \int \D_i^{\ell+1} g(x,\xi,\eta)\big[ \D_i h_{\xi+\alpha e_i}(x) h_{\eta+\beta e_i}(x)-h_{\xi+\alpha e_i}(x)\D_i h_{\eta+\beta e_i}(x)\big]dx \\
&\hspace{120pt} + 2(\beta-\alpha) \int \D_i^\ell g(x,\xi,\eta)h_{\xi+\alpha e_i}(x) h_{\eta+\beta e_i}(x)dx \Big\} \\
&\hspace{40pt}=\sum_{\alpha,\beta,\ell,N} C(N,\xi_i,\eta_i,\alpha,\beta) \\
&\hspace{80pt}\times\Big\{ -a_{\xi_i+\alpha+1} \int \D_i^{\ell+1} g(x,\xi,\eta) h_{\xi+(\alpha+1)e_i}(x)h_{\eta+\beta e_i}(x)dx \\
&\hspace{100pt} +a_{\xi_i+\alpha} \int \D_i^{\ell+1} g(x,\xi,\eta) h_{\xi+(\alpha-1)e_i}(x)h_{\eta+\beta e_i}(x)dx \\
&\hspace{100pt} +a_{\eta_i+\beta+1} \int \D_i^{\ell+1} g(x,\xi,\eta) h_{\xi+\alpha e_i}(x)h_{\eta+(\beta+1) e_i}(x)dx \\
&\hspace{100pt} -a_{\eta_i+\beta} \int \D_i^{\ell+1} g(x,\xi,\eta) h_{\xi+\alpha e_i}(x)h_{\eta+(\beta-1) e_i}(x)dx \\
&\hspace{100pt} +2(\beta-\alpha)\int \D_i^{\ell} g(x,\xi,\eta) h_{\xi+\alpha e_i}(x)h_{\eta+\beta e_i}(x)dx \Big\}
\end{align*}
where the sum $\sum_{\alpha,\beta,\ell,N}$ runs over 
\begin{align*}
-N\le \alpha,\beta\le N,&& 1\le\ell\le N, && |\alpha-\beta|\le N.
\end{align*}
Now since $\alpha,\beta\le N$ and $|\alpha-\beta|\le N$ we see that for all $\mu\in\{\xi_i+\alpha+1,\xi+\alpha,\eta_i+\beta+1,\eta_i+\beta\}$, 
$$ |a_\mu|\lesi \ip{\xi_i\vee \eta_i}^{1/2}, $$ 
and so by the inductive hypthesis for $C(N,\xi_i,\eta_i,\alpha,\beta)$ we have
\begin{align*}
|a_\mu| \big|C(N,\xi_i,\eta_i,\alpha,\beta)\big|  
\lesi \ip{\xi_i\vee \eta_i}^{1/2} \big|C(N,\xi_i,\eta_i,\alpha,\beta)\big|
\lesi \ip{\xi_i\vee \eta_i}^{(N+1)/2}.
\end{align*}
We also have
$$ |2(\beta-\alpha)|\big|C(N,\xi_i,\eta_i,\alpha,\beta)\big| \lesi \ip{\xi_i\vee\eta_i}^{N/2}\le \ip{\xi_i\vee\eta_i}^{(N+1)/2}.$$
This yields \eqref{eq:IBP2}.
For $(\wt{\alpha},\wt{\beta})\in\{(\alpha+1,0), (\alpha-1,0), (\alpha,\beta+1), (\alpha,\beta-1), (\alpha,\beta)\}$ we also have
$$ |\wt{\alpha}-\wt{\beta}|\le |\alpha-\beta|+1\le N+1.$$
Furthermore, the new sum (indexed over $\wt{\alpha}$ and $\wt{\beta}$) for $N+1$ contains 5 new terms for each term in the sum for $N$, and hence we have at most $5^{N+1}$ terms in total. Thus we have obtained the required result for the $N+1$ case which concludes our proof of Proposition \ref{prop:IBP}.

\subsection{Decomposition of symbols and kernel estimates}\label{sec:decomp}
In this section we describe a Littlewood--Paley type decomposition for pseudo-multipliers and the resulting estimates on their kernels. This decomposition will be needed in the proofs of our main results. 

Given any real-valued  bump function $\vph \in C^\infty(\RR_+)$ we set $\vph_j(\xi):=\vph\big(2^{-j}\sqrt{2|\xi|+n}\big)$ for $\xi \in \NN_0^n$. Since the Hermite functions $h_\xi$ with $\xi\in \NN_0^n$ are members of $\sz(\RR^n)$ then we may define the operators $\vph_j(\sqrt{\LL})$ on $\sz'(\RR^n)$ by
$$ \vph_j(\sqrt{\LL})f(x)=\sum_{\xi\in\NN_0^n}\vph_j(\xi)  \ip{f,h_\xi}h_\xi(x)\qquad \forall f\in\sz'(\RR^n), x\in \RR^n,$$
where  $ \ip{f,\phi}=f(\phi)$ for $f\in\sz'(\RR^n)$. 

Throughout the rest of the article we shall fix a bump function $\vph \in C^\infty(\RR_+)$ such that
\begin{align}\label{eq:admissible1}
\supp\vph\subset \big[\tfrac{1}{4},1\big], \qquad\text{and}
\end{align}
\begin{align}\label{eq:admissible2} \sum_{j\ge 0} \vph(2^{-j}\lambda)=1,\qquad \lambda \ge \tfrac{1}{2}.\end{align}

In view of the support of $\varphi_j,$ it follows that 
$$ \vph_j(\sqrt{\LL})f(x)=\sum_{\xi\in \I_j}\vph_j(\xi) \ip{f,h_\xi}h_\xi(x),$$
where
$$\I_j=\big\{\xi\in\NN_0^n: \tfrac{1}{2}4^{j-2}-\tfrac{n}{2} \le |\xi|\le  \tfrac{1}{2}4^{j}-\tfrac{n}{2}\big\},\qquad j\in \NN_0.$$

In addition we can express any symbol $\sigma:\RR^n\times\NN_0^n\to \mathbb{C}$ through the resolution
\begin{align*}
\sigma(x,\xi) = \sum_{j\ge 0}\sigma(x,\xi)\vph_j(\xi),\qquad (x,\xi)\in \RR^n\times\NN_0^n.
\end{align*}
Accordingly, one then has the decomposition for the pseudo-multiplier  \eqref{eq:hpdo}:
\begin{align}\label{eq:LPdecomp}\sigma(\cdot,\LL) =\sum_{j\ge 0}\sigma(\cdot,\LL) \vph_j(\sqrt{\LL}) =\sum_{j\ge 0}\sigma_j(\cdot,\LL)\end{align}
where $\sigma_j(\cdot,\LL)$ is the pseudo-multiplier with symbol $\sigma_j(x,\xi)=\sigma(x,\xi)\vph_j(\xi)$ given by 
\begin{align*}
\sigma_j(x,\LL)f(x)=\sum_{\xi\in\I_j} \sigma_j(x,\xi)\ip{f,h_\xi}h_\xi(x) = \int K_j(x,y)f(y)\,dy
 \end{align*} 
with kernel
\begin{align}\label{eq:Kj} K_j(x,y)=\sum_{\xi\in \I_j}\sigma_j(x,\xi)h_\xi(x)h_\xi(y).\end{align}
Moreover, note that the adjoint $\sigma_j(\cdot,\LL)^*$ of $\sigma_j(\cdot,\LL)$ is given by 
 \begin{align*}
\sigma_j(x,\LL)^*f(x)=\sum_{\xi\in\I_j}h_\xi(x) \ip{\sigma_j(\cdot,\xi)h_\xi, f}= \int K_j^*(x,y)f(y)\,dy
 \end{align*}
 where
 \begin{align} \label{eq:Kjadjoint} K_j^*(x,y)=\sum_{\xi\in \I_j}\overline{\sigma_j(y,\xi)}h_\xi(x)h_\xi(y).\end{align}

We have the following estimates on these kernels with bandlimited symbols. These estimates will be needed in the proof of Theorem \ref{thm:main1}.

\begin{Lemma}\label{lem:Kj}
Let $\K\in\NN_0$ and $m\in\RR$. Assume that the symbol $\sigma: \RR^n\times \NN_0^n\to\mathbb{C}$ satisfies
\begin{align*}
\big\Vert\diff_\xi^\kappa \sigma(\cdot,\xi)\big\Vert_{L^\infty(\RR^n)}\le C_\kappa \ip{\xi}^{m-|\kappa|},\qquad |\kappa|\le \K, \;\xi\in\NN_0^n.
\end{align*}
Then for each $0\le M\le \K$ there exists $C>0$ such that
\begin{align*}
\sup_x\big\Vert (x-\cdot)^\gamma K_j(x,\cdot) \big\Vert_{L^2}\le C\, 2^{j(m+\f{n}{2}-M)},\qquad |\gamma|=M.
\end{align*}
\end{Lemma}

\noindent The proof of Lemma \ref{lem:Kj} requires the following.
\begin{Lemma}[Lemma 2.2 in \cite{LN21}]\label{lem: hoppe}
Let $\phi$ be a smooth function defined in $[0,\infty),$ set $\phi_j(x)=\phi(2^{-j}x)$ for $j\in \NN_0$ and consider  $\ell, N\in \NN,$ $N\ge \ell.$
If $\phi^{(m)}(0)=0$ for all $m\in \NN,$ it holds that 
\begin{equation*}
 |\diff^\ell_\xi(\phi_j(\sqrt{2|\xi|+n}))|\lesssim   \|\phi^{(N)}\|_{L^\infty}  2^{-jN} \ip{\xi}^{N/2-\ell}\quad \forall j,\xi\in \NN_0,
 \end{equation*}
 where the implicit constant depends on $N$, $n$ and $\ell.$
\end{Lemma}

\begin{proof}[Proof of Lemma \ref{lem:Kj}]

Let us first consider the case $M=0$. Then by orthogonality of the hermite functions,
\begin{align*}
\big\Vert  K_j(x,\cdot) \big\Vert_{L^2}^2
= \Bip{\sum_{\xi\in\I_j} \sigma_j(x,\xi)h_\xi(x)h_\xi(\cdot), \sum_{\eta\in\I_j} \sigma_j(x,\eta)h_\eta(x)h_\eta(\cdot)}
=\sum_{\xi\in\I_j} |\sigma(x,\xi)|^2 \vph_j(\xi)^2 h_\xi(x)^2.
\end{align*}
Our assumption on the symbols $\sigma$ then gives
\begin{align*}
\big\Vert  K_j(x,\cdot) \big\Vert_{L^2}^2
\lesi \sum_{\xi\in \I_j} \ip{\xi}^m \,h_\xi(x)^2 
\lesi 4^{jm}\QQ_{4^j + M}(x,x) 
\lesi 2^{j(2m+n)},
\end{align*}
and taking square roots gives the required estimate.

We now consider the case $M\ge 1$. It suffices to prove 
\begin{align}\label{eq:plancherelproof1}
\int |(x_i - y_i)^M K_j(x,y)|^2\,dy \lesi 2^{j(2m+n-2M)}, \qquad \forall\, 1\le i\le n
\end{align}
with constants independent of $i$. Indeed we then have by convexity of the function $t^M$ for $M\ge 1$, that whenever $|\gamma|=M$,
\begin{align*}
\int |(x-y)^\gamma K_j(x,y)|^2\,dy
&=\sum_{i=1}^n \int \big|(x_i - y_i)^2\big|^M |K_j(x,y)|^2\,dy \\
&\lesi \sum_{i=1}^n \int \big|(x_i - y_i)^{2M}\big| |K_j(x,y)|^2\,dy\\
&= \sum_{i=1}^n\int |(x_i - y_i)^M K_j(x,y)|^2\,dy
\end{align*}
which yields the desired result  in view of \eqref{eq:plancherelproof1}.

Let us  continue with the proof of \eqref{eq:plancherelproof1}. Fix $1\le i\le n$ and $M\ge 1$. Then by Proposition \ref{prop:FDI}
\begin{align*}
&\int |(x_i - y_i)^M K_j(x,y)|^2\,dy \\
&\hspace{60pt}=\Big\langle(x_i - (\cdot)_i)^M K_j(x,\cdot), (x_i - (\cdot)_i)^M K_j(x,\cdot)\Big\rangle \\
&\hspace{60pt}=2^{-2M} \sum_{\xi,\ell,\nu,\om} \sum_{\mu,k,\alpha,\beta} c_{\nu, \ell,M}c_{\alpha,k,M} d_\nu(\xi_i)d_\om(\xi_i) d_\alpha(\mu_i)d_\beta(\mu_i)\\
&\hspace{140pt}\times \diff_i^\ell\sigma_j(x,\xi)\diff_i^k \sigma_j(x,\mu) h_{\xi+\nu e_i}(x)h_{\mu+\alpha e_i}(x) \big\langle h_{\xi+\om e_i} ,h_{\mu+\beta e_i} \big\rangle
\end{align*}
where  the sums run over
\begin{align}\label{eq:plancherelproof2}
\xi\in\I_j, &&
\f{M}{2}\le \ell\le M, &&
\nu,\om\ge 0, &&
\nu +\om =  2\ell-M,
\end{align}
and
\begin{align}\label{eq:plancherelproof3}
\mu\in\I_j, &&
\f{M}{2}\le k\le M, &&
\alpha,\beta\ge0,&&
\alpha+\beta = 2k-M.
\end{align}
By orthogonality of the hermite functions the inner product vanishes
unless 
\begin{align}\label{eq:plancherelproof4}
\xi+\om e_i=\mu+\beta e_i\end{align}
and so the expression reduces to 
\begin{align*}
&\int |(x_i - y_i)^M K_j(x,y)|^2\,dy  \\
&\hspace{60pt}=2^{-2M} \sum_{\substack{\xi,\ell, \nu,\om\\ \mu, k ,\alpha,\beta}}  c_{\nu, \ell,M}c_{\alpha, k} d_\nu(\xi_i)d_\om(\xi_i) d_\alpha(\mu_i)d_\beta(\mu_i) \\
&\hspace{140pt} \times \diff_i^\ell\sigma_j(x,\xi)\diff_i^k \sigma_j(x,\mu) h_{\xi+\nu e_i}(x)h_{\mu+\alpha e_i}(x)
\end{align*}
with the sum running over \eqref{eq:plancherelproof2}, \eqref{eq:plancherelproof3} and \eqref{eq:plancherelproof4}. 
By the Cauchy--Schwarz inequality and the symmetry of the resulting terms, we then have
\begin{align*}
&\int |(x_i - y_i)^M K_j(x,y)|^2\,dx  \\
&\hspace{80pt}\le 2^{-2M}\Big(\sum_{\substack{\xi,\ell, \nu,\om\\ \mu, k ,\alpha,\beta}} \big|c_{\nu,\ell,M} d_\nu(\xi_i)d_\om(\xi_i)\diff_i^\ell \sigma_j(x,\xi) h_{\xi+\nu e_i}(x)\big|^2\Big)^{1/2} \\
&\hspace{140pt}\times \Big(\sum_{\substack{\xi,\ell, \nu,\om\\ \mu, k ,\alpha,\beta}} \big|c_{\alpha,k,M} d_\alpha(\mu_i)d_\beta(\mu_i)\diff_i^k \sigma_j(x,\mu) h_{\mu+\om e_i}(x)\big|^2\Big)^{1/2}  \\
&\hspace{80pt}= 2^{-2M}\sum_{\substack{\xi,\ell, \nu,\om\\ \mu, k ,\alpha,\beta}} \big|c_{\nu,\ell,M} d_\nu(\xi_i)d_\om(\xi_i)\diff_i^\ell \sigma_j(x,\xi) h_{\xi+\nu e_i}(x)\big|^2.
\end{align*}
Next, from the Leibniz formula for finite differences (see \cite[Lemma B.2]{LN21}) we have
\begin{align*}
 \big|\diff_i^\ell \sigma_j(x,\xi)\big|
 =\Big|\diff_i^\ell \big[\sigma(x,\xi)\vph_j(\xi)\big]\Big| 
\le \sum_{r=0}^\ell \tbinom{\ell}{r} \big|\diff_i^r(\vph_j(\xi))\big| \, \big|\diff_i^{\ell-r}\sigma(x,\xi+re_i)\big|.
\end{align*}
Now, invoking Lemma~\ref{lem: hoppe} we have
\begin{align*}
\big|\diff_i^r\vph_j(\xi)\big| \lesssim  2^{-jM}\ip{\xi}^{M/2-r},
\end{align*}
and our assumption on $\sigma$ gives
\begin{align*}
\big|\diff_i^{\ell-r}\sigma(x,\xi+re_i)\big| 
\lesssim  \ip{\xi}^{m/2+r-\ell}.
\end{align*} 
These last three facts gives
\begin{align}\label{eq:FDIproof1}
 \big|\diff_i^\ell \sigma_j(x,\xi)\big|
\lesssim 2^{-jM} \ip{\xi}^{m/2+M/2-\ell} .
\end{align}
By Remark \ref{rem:FDI} we also have 
\begin{align}\label{eq:FDIproof2} |d_\nu(\xi_i)d_\om(\xi_i)|\lesi \ip{\xi}^{(\nu+\om)/2}=\ip{\xi}^{\ell-M/2}.
\end{align}
Then from \eqref{eq:FDIproof1} and \eqref{eq:FDIproof2}  (and bearing in mind $\ip{\xi}\sim 4^j$ for $x\in\I_j$) we have
\begin{align*}
\int |(x_i - y_i)^M K_j(x,y)|^2\,dx 
\lesi 2^{2j(m-M)} \sum_{\substack{\xi,\ell, \nu,\om\\ \mu, k ,\alpha,\beta}}h_{\xi+\nu e_i}(x)^2.
\end{align*}
Note that the summation in $\mu, k,\alpha,\beta$ only add at most $O(M)$ terms to each item in the sum. Indeed from \eqref{eq:plancherelproof2}--\eqref{eq:plancherelproof4} we have 
$$|\xi_i-\mu_i| =|2\om-2\beta|\le 4\ell-2M \le 2M.$$
Thus we have
\begin{align*}
\int |(x_i - y_i)^M K_j(x,y)|^2\,dx 
\lesi 2^{2j(m-M)}\sum_{|\xi|\le 4^j+2M}h_{\xi+\nu e_i}(x)^2
\lesi 2^{2j(m-M+\f{n}{2})}
\end{align*}
in view of \eqref{eq:QQ} and \eqref{QQ bound}.
This completes the proof of \eqref{eq:plancherelproof1}, and hence the lemma.
\end{proof}

\section{$L^2$ boundedness}\label{sec:L2}

\subsection{Proof of Theorem \ref{thm:main1}}

The proof of Theorem \ref{thm:main1} uses the well known Cotlar--Knapp--Stein Lemma (see for example \cite[Lemma 8.5.1]{Graf2}).
\begin{Lemma}[Cotlar--Knapp--Stein]\label{lem:CKS} 
Let $\{T_j\}_{j\in\ZZ}$ be a family of operators satisfying:
\begin{enumerate}[\upshape(i)]
\item Each $T_j$ is bounded on $L^2(\RR^n)$.
\item For each  $j,k\in\ZZ$,
$$ \Vert T_j^* T_k\Vert_{L^2\to L^2} +\Vert T_j T_k^*\Vert_{L^2\to L^2} \le \gamma(j-k)$$
where $\gamma$ satisfies $\sum_{j\in\ZZ}\gamma(j-k)=C<\infty$. 
\end{enumerate}
Then we have 
\begin{enumerate}[\upshape (a)]
\item For all finite subset $\Lambda$ of $\ZZ$ we have $\Vert \sum_{j\in\Lambda} T_j\Vert_{L^2\to L^2}\le C$. 
\item For each $f\in L^2(\RR^n)$ we have $\sum_{j\in\ZZ} \Vert T_j f\Vert_{L^2}^2 \le C^2 \Vert f\Vert_{L^2}^2$.
\item For each $f\in L^2(\RR^n)$ the sequence $\sum_{|j|\le N} T_jf $ converges to some $Tf$ as $N\to\infty$ in the norm topology of $L^2(\RR^n)$. The linear operator defined this way is bounded from $L^2(\RR^n)$ to $L^2(\RR^n)$ with $\Vert T\Vert_{L^2\to L^2}\le C$. 
\end{enumerate}
\end{Lemma}

Assume firstly that $\sigma$ compactly supported in the spatial variable; we will obtain bounds that are independent of the support of $\sigma$ in $x$. Let $\vph$ be a bump function satisfying  \eqref{eq:admissible1}-\eqref{eq:admissible2} and $\sigma_j(\cdot,\LL)$ be defined as in \eqref{eq:LPdecomp} with kernel in \eqref{eq:Kj}. We shall apply the Cotlar--Knapp--Stein lemma to the family $\{\sigma_j(\cdot,\LL)\}_{j\ge 0}$.

We begin by showing condition (i) of Lemma \ref{lem:CKS}:  more precisely we show that there exists $C_0>$ independent of $j$ such that 
\begin{align}\label{eq:RTPC1}
\Vert \sigma_j(\cdot,\LL)\Vert_{L^2\to L^2} \le C_0, \qquad\forall j\ge 0.
\end{align}
Taking $\K=\floor{n/2}+1$ we have, by the Cauchy--Schwarz inequality,
\begin{align*}
\Vert \sigma_j(\cdot,\LL)f\Vert_{L^2}^2 
\le \int \Big(\int | K_j(x,y) |^2(1+2^j|x-y|)^{2\K}dy\Big)\Big(\int \f{|f(y)|^2}{(1+2^j|x-y|)^{2\K}}dy\Big)dx.
\end{align*}
Firstly, applying Lemma \ref{lem:Kj}  with  $M=0$ we obtain
\begin{align*}
\int_{|y-x|\le 2^{-j}}| K_j(x,y) |^2 (1+2^j|x-y|)^{2\K} \,dy
\lesi \int  |K_j(x,y)|^2\,dy 
\lesi 2^{jn}.
\end{align*}
On the other hand, applying Lemma \ref{lem:Kj}  with $M=\K=$ we have
\begin{align*}
\int_{|y-x|>2^{-j}}| K_j(x,y) |^2 (1+2^j|x-y|)^{2\K}\,dy 
\lesi 2^{2\K j} \int |x-y|^{2\K} |K_j(x,y)|^2\,dy 
\lesi 2^{jn}.
\end{align*}
Combining both these estimates along with the fact that $\K>n/2$ we have
\begin{align*}
\Vert \sigma_j(\cdot,\LL)f\Vert_{L^2}^2
\lesi 2^{jn}\int |f(y)|^2\int \f{dx}{(1+2^{j}|x-y|)^{2\K}}\, dy
\lesi  \Vert f\Vert_{L^2}^2, 
\end{align*}
which gives \eqref{eq:RTPC1}.

Secondly, we show that there exists $C_1>0$ and $\ep>0$ such that for each  $j,k\in\ZZ$,
\begin{align}\label{eq:RTPC2}
 \Vert \sigma_j(\cdot,\LL)^* \sigma_k(\cdot,\LL)\Vert_{L^2\to L^2} +\Vert \sigma_j(\cdot,\LL) \sigma_k(\cdot,\LL)^*\Vert_{L^2\to L^2} \le C_1 \,2^{-\ep|j-k|}.
 \end{align}
With both \eqref{eq:RTPC1} and \eqref{eq:RTPC2} in hand we may invoke the Cotlar--Knapp--Stein lemma to conclude the result of our theorem.

 Let us continue with the proof of \eqref{eq:RTPC2}. 
We first consider the family $\big\{\sigma_j(\cdot,\LL)\sigma_k(\cdot,\LL)^*\big\}_{j,k\ge 0}$. From the supports of $\vph_j$ and $\vph_k$ we have $ \I_j\cap \I_k =\emptyset$ whenever $|j-k|\ge 2$. Hence  for $|j-k|\ge 2$ we have
\begin{align}\label{eq:proofC.1} 
\sigma_j(\cdot,\LL)\sigma_k(\cdot,\LL)^* 
= \big(\sigma(\cdot,\LL) \vph_j\big)\big(\sigma(\cdot,\LL) \vph_k\big)^* 
= \sigma(\cdot,\LL) \vph_j \,\vph_k \sigma(\cdot,\LL)^* = 0.
\end{align}
 If $|j-k|\le 1$ then from \eqref{eq:RTPC1} we have
\begin{align}\label{eq:proofC.2} \Vert \sigma_j(\cdot,\LL)\sigma_k(\cdot,\LL)^*\Vert_{L^2}\le \Vert \sigma_j(\cdot,\LL)\Vert_{L^2}\Vert \sigma_k(\cdot,\LL)^*\Vert_{L^2} = \Vert \sigma_j(\cdot,\LL)\Vert_{L^2}\Vert \sigma_k(\cdot,\LL)\Vert_{L^2} \le C_0^2. \end{align}
Both \eqref{eq:proofC.1} and \eqref{eq:proofC.2} together give \eqref{eq:RTPC2} for the operators $\big\{\sigma_j(\cdot,\LL)\sigma_k(\cdot,\LL)^*\big\}_{j,k\ge 0}$.

We turn to the operators $\big\{\sigma_j(\cdot,\LL)^*\sigma_k(\cdot,\LL)\big\}_{j,k\ge 0}$. For $|j-k|\le 1$, as in \eqref{eq:proofC.2} we may use \eqref{eq:RTPC1} to see that
\begin{align}\label{eq:proofC.3} 
\Vert \sigma_j(\cdot,\LL)^*\sigma_k(\cdot,\LL)\Vert_{L^2}
\le \Vert \sigma_j(\cdot,\LL)\Vert_{L^2}\Vert \sigma_k(\cdot,\LL)\Vert_{L^2} 
\le C_0^2. \end{align}
It remains to check the case $|j-k|\ge 2$. Before proceeding further let us consider the expressions for the relevant kernels. By recalling the operators $\sigma_j(\cdot,\LL)$ and $\sigma_j(\cdot,\LL)^*$ and their kernels \eqref{eq:Kj}--\eqref{eq:Kjadjoint} we see that 
 \begin{align*}
\sigma_j(\cdot,\LL)^*\sigma_k(\cdot,\LL)f(x)=\int K_{(j,k)}(x,y)f(y)\,dy
 \end{align*}
 where 
  \begin{align*}K_{(j,k)}(x,y)=\sum_{\eta\in \I_j}\sum_{\xi\in\I_k}h_\xi(x)h_\eta(y)\int \sigma_k(z,\xi)\overline{\sigma_j(z,\eta)}h_\xi(z)h_\eta(z)\,dz.
  \end{align*}
We now continue with the proof of \eqref{eq:proofC.2}.
By orthogonality of the hermite functions we have
\begin{align*}
\sigma_j(x,\LL)h_\eta(x) = \sigma_j(x,\eta)h_\eta(x)
\end{align*}
and thus we see that
\begin{align*}
\ip{\sigma_j(\cdot,\LL)^* \sigma_k(\cdot,\LL) f,h_\eta} 
= \ip{\sigma_k(\cdot,\LL)f, \sigma_j(\cdot,\LL) h_\eta}
= \sum_{\xi\in\I_k} \ip{f,h_\xi} \bip{\sigma_k(\cdot,\xi)h_\xi,\sigma_j(\cdot,\eta)h_\eta}
\end{align*}
Then by Parseval's identity and the Cauchy--Schwarz inequality,
\begin{align}\label{eq:proofC.4}
\Vert \sigma_j(\cdot,\LL)^*\sigma_k(\cdot,\LL) \Vert_{L^2}^2
&= \Vert \ip{\sigma_j(\cdot,\LL)^*\sigma_k(\cdot,\LL) f,h_\eta}\Vert_{\ell^2(\eta)}^2 \\\notag
&\le \sum_{\eta\in\I_j}\Big(\sum_{\xi\in\I_k} |\ip{f,h_\xi}|^2\Big)\Big(\sum_{\xi\in\I_k} \big|\bip{\sigma_k(\cdot,\xi)h_\xi,\sigma_j(\cdot,\eta)h_\eta}\big|^2\Big)
\end{align}
Now we claim that  the following estimate holds for $N=\floor{2n/(1-\dd)}+1$:
\begin{align}\label{eq:proofC.AO}
\big|\bip{\sigma_k(\cdot,\xi)h_\xi,\sigma_j(\cdot,\eta)h_\eta}\big| \lesi \f{\bip{|\xi|\vee|\eta|}^{\f{N}{2}(1+\dd)}}{|\xi-\eta|^N},\qquad \xi\ne\eta.
\end{align}
Assume the claim for the moment. 
Observe that $|j-k|\ge 2$ implies that $|\xi-\eta|\sim 4^{j\vee k}$ for $\eta\in\I_j$ and $\xi\in\I_k$. Observe also that $|\I_j|\sim 4^{jn}$. These facts in conjunction with \eqref{eq:proofC.AO} give
\begin{align*}
\sum_{\eta\in\I_j}\sum_{\xi\in\I_k} \big|\bip{\sigma_k(\cdot,\xi)h_\xi,\sigma_j(\cdot,\eta)h_\eta}\big|^2 
&\lesi \sum_{\eta\in\I_j}\sum_{\xi\in\I_k} \bip{|\xi|\vee|\eta|}^{N(1+\dd)}|\xi-\eta|^{-2N} \\
&\sim \sum_{\eta\in\I_j}\sum_{\xi\in\I_k} 4^{(j\vee k)(N(\dd-1))} \\
&\sim 4^{(j\vee k)(N(\dd-1)+2n)}
\end{align*}
Now set $\ep:=N(1-\dd)-2n>0$.
On inserting the previous estimate into \eqref{eq:proofC.4} we arrive at
\begin{align*}
\Vert \sigma_j(\cdot,\LL)^*\sigma_k(\cdot,\LL) \Vert_{L^2}^2
\le \Big(\sum_{\eta\in\I_j}\sum_{\xi\in\I_k} \big|\bip{\sigma_k(\cdot,\xi)h_\xi,\sigma_j(\cdot,\eta)h_\eta}\big|^2\Big)\sum_{\xi\in\I_k} |\ip{f,h_\xi}|^2 
\lesi  2^{-\ep|j-k|}\sum_{\xi\in\I_k} |\ip{f,h_\xi}|^2
\end{align*}
which gives \eqref{eq:RTPC2} for $\big\{\sigma_j(\cdot,\LL)^*\sigma_k(\cdot,\LL)\big\}_{|j-k|\ge 2}$, modulo claim \eqref{eq:proofC.AO}. 

 To complete the proof it remains to  show the claim \eqref{eq:proofC.AO}. 
Let $i\in\{1,\dots,n\}$. By Proposition \ref{prop:IBP}, our assumption on $\sigma$, and the Cauchy--Schwarz inequality,
\begin{align*}
&|\xi_i-\eta_i|^N\Big|\int \sigma_k(x,\xi)\overline{\sigma_j(x,\eta)}h_\xi(x)h_\eta(x)\,dx\Big|\\
&\qquad\qquad\sim \Big|\sum_{\substack{-N\le \alpha,\beta\le N,\\ 1\le \ell\le N,\\ |\alpha-\beta|\le N}} C(N,\xi_i,\eta_i,\alpha,\beta) \int \D_i^\ell\big(\sigma_k(x,\xi)\overline{\sigma_j(x,\eta)}\big) h_{\xi+\alpha e_i}(x) h_{\eta+\beta e_i}(x)\,dx\Big|.
\end{align*}
By the Leibniz rule, our assumption on the symbol $\sigma$ and the fact that $\ell\le N$, we have
\begin{align*}
\big|\D_i^\ell\big(\sigma_k(x,\xi)\overline{\sigma_j(x,\eta)}\big)\big|
&\le \sum_{s+t=\ell} \tbinom{\ell}{s,t} \big|\D_i^s \sigma(x,\xi)\D_i^t\overline{\sigma(x,\eta)}\vph_k(\xi)\vph_j(\eta)\big| \\
&\lesi \sum_{s+t=\ell} \ip{\xi_i}^{\dd s/2} \ip{\eta_i}^{\dd t/2} \\
&\lesi \ip{\xi_i\vee\eta_i}^{\dd \f{N}{2}}.
\end{align*}
Inserting this estimate into our previous expression and recalling the bounds on $C(N,\xi_i,\eta_i,\alpha,\beta)$ we obtain
\begin{align*}
&|\xi_i-\eta_i|^N\Big|\int \sigma_k(x,\xi)\overline{\sigma_j(x,\eta)}h_\xi(x)h_\eta(x)\,dx\Big| \\
&\qquad\qquad\qquad\lesi \ip{\xi_i\vee\eta_i}^{\f{N}{2}(1+\dd)} \sum_{\substack{-N\le \alpha,\beta\le N,\\ 1\le \ell\le N,\\ |\alpha-\beta|\le N}} \Vert h_{\xi+\alpha e_i}\Vert_{L^2}^2\Vert h_{\eta+\beta e_i}\Vert_{L^2}^2\\
&\qquad\qquad\qquad\lesi \ip{\xi_i\vee\eta_i}^{\f{N}{2}(1+\dd)} 
\end{align*}
by the Cauchy--Schwarz inequality. Then summing over all $i\in\{1,\dots,n\}$ we have
\begin{align*}
&|\xi-\eta|^N\Big|\int \sigma_k(x,\xi)\overline{\sigma_j(x,\eta)}h_\xi(x)h_\eta(x)\,dx\Big| \\
&\qquad\qquad\qquad\lesi
\sum_{1\le i\le n}|\xi_i-\eta_i|^N\Big|\int \sigma_k(x,\xi)\overline{\sigma_j(x,\eta)}h_\xi(x)h_\eta(x)\,dx\Big| \\
&\qquad\qquad\qquad\lesi \sum_{1\le i\le n}\ip{\xi_i\vee\eta_i}^{\f{N}{2}(1+\dd)} \\
&\qquad\qquad\qquad\lesi \bip{|\xi|\vee|\eta|}^{\f{N}{2}(1+\dd)}
\end{align*}
which gives \eqref{eq:proofC.AO}.  This completes the proof Theorem \ref{thm:main1}.

\subsection{Proof of Theorem \ref{thm:main2}}
We shall employ again the decomposition  $\sigma_j(\cdot,\LL)=\sigma(\cdot,\LL)\vph_j(\LL)$ as in \eqref{eq:LPdecomp}.
Our aim is to show that there exists $C>$ independent of $j$ such that 
\begin{align}\label{eq:RTPC3}
\Vert \sigma_j(\cdot,\LL)\Vert_{L^2\to L^2} \le C 2^{jm}, \qquad\forall j\ge 0.
\end{align}
Since $m<0$ then applying \eqref{eq:RTPC3} we have
$$ \Vert \sigma(\cdot,\LL) f\Vert_{L^2} \le \sum_{j\ge 0}\Vert \sigma_j(\cdot,\LL) f\Vert_{L^2} \lesi \sum_{j\ge 0} 2^{jm} \Vert f\Vert_{L^2}\lesi \Vert f\Vert_{L^2}.$$
Let us now prove \eqref{eq:RTPC3}, which is a straightforward extension of the argument used for \eqref{eq:RTPC1}. In fact, applying Lemma \ref{lem:Kj} (firstly with  $M=0$ and secondly with $M=\K:=\floor{n/2}+1$) we obtain
\begin{align*}
\int |K_j(x,y) |^2 (1+2^j|x-y|)^{2\K}\,dy 
\lesi 2^{j(2m+n)}.
\end{align*}
With this estimate in hand we can invoke the Cauchy--Schwarz inequality to see that
\begin{align*}
\Vert \sigma_j(\cdot,\LL)f\Vert_{L^2}^2 
\lesi 2^{j(2m+n)} \int |f(y)|^2\int \f{dx}{(1+2^{j}|x-y|)^{2\K}}\, dy
\lesi 2^{2jm}\Vert f\Vert_{L^2}^2,
\end{align*}
yielding \eqref{eq:RTPC3}.

\subsection{Proof of Theorem \ref{thm:main3}}

In this section we give the proof of Theorem \ref{thm:main3}. 
The approach here is inspired from the toroidal case (see \cite[Theorem 4.8.1]{RT}). 

First set
$$Sf(y,x):= \sum_{\xi\in\NN_0^n} \sigma(y,\xi) h_\xi(x) \ip{f,h_\xi}.$$
Then clearly $\sigma(x,\LL) f(x)=Sf(x,x)$, and that 
\begin{align*}
\Vert \sigma(\cdot,\LL) f\Vert_{L^2}^2 
=\Vert Sf(\cdot,\cdot)\Vert_{L^2}^2
\le \int \Vert Sf(\cdot,x)\Vert_{L^\infty}^2\,dx.
\end{align*}
Let $\N=\floor{n/2}+1$. Then by the well known Sobolev embedding (\cite[Theorem 6.2.4]{Graf2}), it follows that
\begin{align*}
\Vert \sigma(\cdot,\LL) f\Vert_{L^2}^2 
\lesi \int \sum_{|\nu|\le \,\N} \Vert \partial^\nu Sf(\cdot,x)\Vert_{L^2}^2\,dx.
\end{align*}
Now from Fubini's theorem and Parseval's identity we have
\begin{align*}
\int \sum_{|\nu|\le\N} \Vert \partial^\nu_ySf(\cdot,x)\Vert_{L^2}^2\,dx
=\sum_{|\nu|\le\N} \int \Vert \partial^\nu_ySf(y,\cdot)\Vert_{L^2}^2 dy
=\sum_{|\nu|\le\N} \int \big\Vert \bip{\partial^\nu_ySf(y,\cdot),h_\xi}\big\Vert_{\ell^2(\xi)}^2dy.
\end{align*}
Orthogonality of the Hermite functions then gives
\begin{align*}
\big|\bip{\partial^\nu_y Sf(y,\cdot),h_\xi} \big|
\le\sum_{\eta\in\NN_0^n} \big|\partial^\nu_y \sigma(y,\eta)\big|\big|\ip{f,h_\eta}\big|\big| \ip{h_\xi,h_\eta}\big|
=\big|\partial^\nu_y \sigma(y,\xi)\big|\big| \ip{f,h_\xi}\big|
\end{align*}
for $(y,\xi)\in\RR^n\times \NN_0^n$. This fact along with Fubini's theorem gives
\begin{align*}
\Vert \sigma(\cdot,\LL) f\Vert_{L^2}^2 
\lesi \sum_{|\nu|\le\N} \int \sum_{\xi\in \NN_0^n} \big|\partial^\nu_y \sigma(y,\xi) \big|^2 \big|\ip{f,h_\xi}\big|^2 dy
= \sum_{\xi\in\NN_0^n} \big|\ip{f,h_\xi}\big|^2 \Big(\sum_{|\nu|\le\N} \big\Vert \partial^\nu\sigma(\cdot,\xi)\big\Vert_{L^2}^2 \Big).
\end{align*}
Our condition \eqref{eq:sobolev-cond} and Parseval's identity now allows us to conclude that
\begin{align*} 
\Vert \sigma(\cdot,\LL) f\Vert_{L^2}^2 
\lesi \sum_{\xi\in\NN_0^n} \big|\ip{f,h_\xi}\big|^2
=\Vert f\Vert_{L^2}^2,
\end{align*}
which completes the proof of the Theorem.

\end{document}